\documentclass{article}
\usepackage[utf8]{inputenc}

\usepackage{authblk, graphicx,tikz-cd, cleveref}
\usepackage{graphicx}
\usepackage{subcaption}
\usepackage{booktabs}
\usepackage{multirow}

\title{High-structure and its impact on performance gaps in large-enrollment frosh courses}
\author[1]{Pedro Morales-Almazan}
\author[1]{Hanjue Pan}
\author[1]{Jasmine Tom}
\affil[1]{University of California, Santa Cruz}
\date{}
\begin{document}
\maketitle
\begin{abstract}
We introduce the concept of structure-regions and their impact on student experience in large-enrollment frosh courses. These regions involve three major components: active learning, low-stakes assignments, and formative assessment. We analyze performance data and show that high-structure courses have a positive impact on class performance by decreasing the performance gap throughout the duration of a course. 
\end{abstract}

\section{Introduction}\label{sec:intro}
Traditionally, large-enrollment courses suffer from many philosophical and logistical hurdles when it comes to teaching. From engagement and communication, to access and scalability, large-enrollment courses face challenges that can be mitigated by integrating a robust structure.

Despite these being mainly lower division courses where students require more guidance, large-enrollment courses tend to be traditional lecture-based with minimal instructor-student interaction \cite{borda2017adapting,geske1992overcoming,talbot2015transforming}. This passive approach is in general ineffective in fostering deep understanding of scientific concepts \cite{borda2017adapting,talbot2015transforming}.
The large number of students, often in the hundreds, makes it challenging to implement research-based teaching practices \cite{talbot2015transforming, borda2017adapting, henderson2007barriers, allen2005infusing, geske1992overcoming}. From the use of active learning techniques to grading and providing meaningful feedback, large-enrollment courses suffer from issues that are non-scalable \cite{allen2005infusing, geske1992overcoming}. 
Student engagement is also difficult in a large-enrollment setting. Due to the physical layout of these classrooms, students perceive this experience to be highly impersonal and intimidating \cite{geske1992overcoming}. Usually these lectures are held in lecture halls where seating arrangements are fixed. This prevents students from interacting with each other and with the instructor, making them inaccessible for meaningful conceptual exchanges \cite{henderson2007barriers, geske1992overcoming}. This reinforces the passive role of students which prevents them from developing deep conceptual understanding \cite{borda2017adapting}.
Communication often feels impersonal and intimidating for students in and outside the classroom \cite{geske1992overcoming}. Inside the classroom, due to the volume and physical layout in lecture halls, instructors tend to be seen as remote figures and students lean to adopt an spectator role \cite{geske1992overcoming}. It is challenging to effectively manage student communication due to the limited human, time, and technological resources \cite{geske1992overcoming}. Outside the classroom, providing meaningful feedback for student work also becomes a difficult challenge, often with a reduced frequency and quality of it \cite{talbot2015transforming}. These characteristics tend to develop poor student motivation and satisfaction \cite{talbot2015transforming}.
The use of technology has helped addressing some of the issues present in large-enrollment science courses \cite{cooper2018influence, mestre2002effect, beichner2007student, hunter2000use}. This provides a way of addressing the issues in scalability, engagement, and feedback in a more efficient way. In this paper we propose tackling these persistent challenges in large-enrollment courses with the framework of \emph{structure}. By this we specifically consider three components: low-stakes assignments, active-learning activities, and formative assessments. By purposefully addressing these three components, a large-enrollment course can be moved from a low-structure region to a high-structure region, improving the student experience.

This article is organized in six sections. \Cref{sec:intro} presents an introduction to the topic including some of the most important challenges of large-enrollment courses. \Cref{sec:structure} describes what we define as structure in a classroom and the concept of structure regions. In \Cref{sec:imp} we describe the implementation of a highly-structured precalculus course at the University of California, Santa Cruz. \Cref{sec:met} describes the analysis made on the impact of high-structured approaches in student performance from two large-enrollment math courses at the University of California, Santa Cruz. \Cref{sec:res} presents the results of the analysis on these two courses. In \Cref{sec:con} we make some final remarks and conclusions. 

\section{Structure}\label{sec:structure}

We define structure as the logistical breakdown of the activities in a course. This falls in line with the notion that highly-structured courses have a higher component of formative assessments and active learning, which has shown to impact student performance \cite{freeman2011increased}. 
    
Structure then can be considered as a spectrum going from low to high. Low-structured courses are those where instruction is held in a traditional lecture style with students playing a passive role, few high-stakes assignments, and no formative assessments. In contrast, high-structured courses include active learning activities, many low-stakes assignments, and multiple formative assessments. 
    
\begin{figure}[h]
    \centering
    \begin{tikzpicture}[commutative diagrams/every diagram,
        declare function={R=4;Rs=R*cos(60);}]
        \path 
        (0,0)  node (D)  {.} 
        (0:0.4*R) node {{\small\emph{low structure}}}
        (135:0.75*R) node {{\small\emph{higher structure}}}
        (-30:R) node[right] {Low-stakes assignments} 
        (210:R) node[left] {Active learning}
        (90:R) node[above] {Formative assessments}
        (-30:R) node (X) {} 
        (90:R) node (Y) {}
        (210:R) node (Z) {};
        \draw[dashed] (0,0) circle (0.8*R);
        \path[commutative diagrams/.cd, every arrow, every label]
        (D) foreach \X in {X,Y,Z} {edge (\X)}; 
    \end{tikzpicture}
    \caption{Dimensions of course structure.}
    \label{fig:strucdim}
\end{figure}
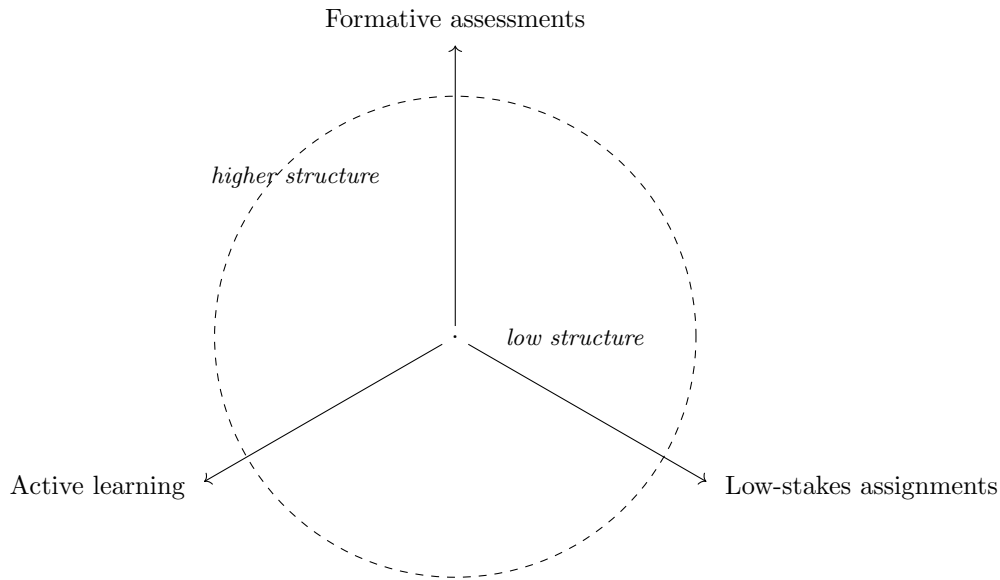

These three components provide a way to identify the level of structure a course have. We can see in \Cref{fig:strucdim} that we can move from a region of low structure to a higher structure by increasing in any of these directions. 
    
It has been shown that high-structure positively impacts student performance with respect to low and medium-structured courses \cite{freeman2011increased, haak2011increased}. Even more, high-structure decreases the equity gap experienced by underrepresented students \cite{haak2011increased}. Other studies suggest the importance on structure even in remote settings as a way to compensate for highly interactive environments \cite{lee2009influence}. This is particularly relevant, not only in distance learning, but also in large-enrollment courses. 
    
We could say that large-enrollment courses behave similarly to distance courses in the sense that students do not experience the same level of interaction and have an \emph{in-mass} experience \cite{cuseo2007empirical, macgregor2000strategies}. Typically, these courses have low engagement and interaction among students and with instructors \cite{cotner2008rapid, swap2015approach, cuseo2007empirical, macgregor2000strategies}.
    
Traditional large-enrollment intro level courses then exist in the \emph{low-structure region}, being heavily lectured, impersonal and with low student interaction, relying on few high-stakes assessments \cite{cuseo2007empirical, macgregor2000strategies, oliver2007exploring, biggs2003enriching}.
    
    \subsection{Impact of structure}
    The amount of active learning activities, formative assessment, and low-stakes assignments relate to the engagement and performance of students in large-enrollment courses. Moving from a region of low-structure into a higher-structure region positively impacts students.
    
    Using the framework of structure it is possible to better identify effective ways of addressing some of the typical issues present in teaching large-enrollment courses. We can organize these issues in three main components: 
    \begin{itemize}
        \item Scalability
        \item Engagement
        \item Feedback
    \end{itemize}
        \subsubsection{Scalability}
        Scalability in general refers to the characteristic of a system or process to adapt to an increasing number of elements and to be capable of enlargement \cite{bondi2000characteristics}. In the case of the teaching, scalability refers to the quality of being able to adapt to an increase in students for a course without diminishing the teaching and learning quality. Both logistical and physical constrains can affect the scalability of courses which in turn affect the student experience. Modern technology and the implementation of teaching teams has made possible to adopt scalable solutions in and outside the classroom. 
    
        A common artifact of non-scalable frameworks is to have few high-stakes assignments in a course. Due to lack of human and logistical resources, instructors often lean towards concentrating their course work and grading distribution on few high-stakes assignments, such as a midterm and a final, and a few of written homework assignments. High-stakes assignments often put a lot of pressure on students. Motivation decreases in high-stakes assessments, specially for minoritized students \cite{amrein2003effects, doi:10.1207/s15430421tip42013, madaus2001adverse}. These also lead to lower retention and higher dropout rates \cite{amrein2003effects}. 
    
        In contrast, low-stakes assignments serve as an effective way to provide content practice and more effective approach in assigning grades. The use of low-stakes assignments provides students with the opportunity to routinely practice and increase their exposure to content. Using multiple low-stakes assignments fosters a more equitable system for assigning grades. Students' grades are taken to be the average of multiple assignments instead of just a single one. This avoids the uncertainty of a few assignments accounting for the student's grade. When considering multiple assignments, the validity of the grading system can be improved since the random errors in grades will partially cancel each other out \cite{hopkins1989errors}.
    
        \subsubsection{Engagement}
        There are many conceptualizations of what engagement means in the classroom setting. Some of the most common ones focus on the level of interaction, participation, and interest to the subject matter from students \cite{ahlfeldt2005measurement,bryson2007role}. Engagement has profound impact in student learning and the development of self-regulated learners \cite{corno1983role}. 
        
        Usually in large-enrollment settings there is very little room for student interaction, participation, and the development of motivation on the subject matter. Student interaction becomes challenging, even more when the physical constrains of the room do not allow for chairs to be moved or rearranged. Also, it is difficult to have students actively participate, defaulting to a more passive role. This also makes it more challenging to develop student motivation towards the topics in the course. 
        
        The incorporation of active learning activities can positively impact engagement and sense of belonging for students in large-enrollment courses. Active learning activities aim for students to play an \emph{active} role in the classroom. This fosters student-student and student-instructor interactions, creating a deeper sense of belonging, more class engagement, and learning communities. 
    
        \subsubsection{Feedback}
        Feedback establishes a closed-loop in the learning process. Effective feedback can be used not only as a way to signaling students on the correctness of their work, but also as a way to guide their cognitive processes \cite{bangert1991instructional, vollmeyer2005surprising}. The expectation of having feedback can alter student strategies towards assignments, helping them acquire more knowledge in fewer attempts \cite{vollmeyer2005surprising}. 
        
        Providing effective feedback becomes challenging in large-enrollment courses. Effective feedback usually requires meta-cognitive assessment and guidance, which can be time consuming and logistically laborious for large groups. Even more for paper-based assignments, returning them to students becomes challenging and sometimes even impossible. 
        
        Formative assessment can be utilized in large-enrollment courses for providing effective feedback and promote non-academic skills that are fundamental for student success. One of the main goals of formative assessment is to monitor the progress of the learning process in students \cite{black1998assessment}. This can be used for providing effective feedback to students that can them to focus on their weaknesses and strengths, as well as identifying learning strategies to enhance their achievement \cite{mccarthy2017enhancing}. Formative assessment develops meta-cognition and self-efficacy in students which are fundamental, not only for the success in a particular course, but in their academic and professional career. 
        
\section{Implementation: Large-enrollment precalculus course}\label{sec:imp}

Precalculus courses at the University of California, Santa Cruz are often large-enrollment: usually between 200-300 students. The setup is traditionally lecture based on stadium-style rooms with optional smaller discussion sections held by teaching assistants. 

One of the authors had been in charge of teaching these precalculus courses in the past few years. Following the high-structure approach, the setup of the course incorporates more active learning components, low-stakes assignments, and more formative assessments. This has been made possible with the implementation of technology and the support of a teaching team involving teaching assistants, learning assistants, and tutors.

    \subsection{Active Learning}
        Implementing active learning activities in any classroom can be challenging at first, even more when the number of students is in the hundreds and the seating is fixed. However, it is possible to adapt to these circumstances and to provide a more active experience for students.
        
        For this course, the instructor includes small discussions to motivate or reflect upon the concepts seen in class. A sample prompt would be ``\emph{write down an example of something that is not a function and talk to your neighbor about it.}'' The instructor would give students a couple of minutes to write and discuss these, and then proceed to follow the main discussion about it. This might happen a few times during a regular lecture, where students are invited to think and share about concepts.
        
        Another activity used is to work problems in groups. This might be challenging in rooms where the seats are bolted on the floor. However, it is still possible to have students work with their adjacent neighbors. That way groups of 2 or 3 students are possible. The instructor would display a problem and give 10-15 minutes for students to work on it. During this time, the teaching team would circulate in the room checking students' work and answering questions. Depending on the seating arrangement, physically accessing students to address their questions might be difficult. When possible, blocking students from seating every other row can give some walking room for the teaching team to access students. 
        
        During discussion section, the main approach from teaching assistants is to cover material using active learning techniques. Discussion sections are smaller sections with 30-35 students held in smaller rooms with better seating arrangements. By making these sessions to be mandatory and by having teaching assistants to focus on different active learning strategies, students experience more practice of the topics covered in class.
        
    \subsection{Formative assessments}
        Fluid two-way communication can be challenging in large-enrollment courses. Giving any type of meaningful feedback becomes a very difficult task when the number of students are in the hundreds. With the help of technology, helping students monitor their own learning progress and giving more effective feedback becomes a more realistic task. 
        
        Most of the communication for this precalculus course is handled through the University's Learning Management System (LMS). This LMS is used to have daily reflections regarding readings that students have to do before each lecture. These reflections are used as a metacognitive device to prompt students to think about the challenging and interesting concepts that will be covered in lecture. 
        
        Also, these reflections provide information to the instructor. Using the data generated by these reflections, the instructor would survey the most common challenges and motivators for each lecture. With this the instructor is able to adapt the lesson plan and give feedback to students in an indirect way.
        
        Another formative assessment tool used for this course is the implementation of a weekly application worksheet during lecture. This is worked using a scaffolding approach. Students are encouraged to work in groups. The worksheet is worked one part at a time, taking around 5-10 min per part. While students work on each part, the teaching team circulates in the room, checking for student progress and answering questions. The instructor then asks the students for ideas on how to solve each part and explains the most important points of the solution. Then moves to the next part. 
    
    \subsection{Low-stakes assignments}
        With a big number of students, managing and grading multiple assignments becomes logistically complex. Even with large teaching teams, grading could be a very costly task. To a certain extent, it is possible to think that there are better ways to spend highly trained resources that could be of more benefit to the students. The use of automatic grading platforms can help with the more efficient usage of resources while giving students the opportunity of more practice.
        
        These platforms are well suited for low-stakes assignments. It is possible to assess computational exercises, multiple-choice questions, and numerical questions while giving accurate and effective feedback to students. 
        
        This feature enabled the incorporation of small exercises in readings that students completed before each lecture day. This gives them the opportunity assess their readings immediately. At the same time, having more of these assignments make each of them to have less weight in their overall grade. This takes pressure off from students and helps them to try more honestly, without focusing on getting the question right, but more on assessing their own understanding. 
        
        Likewise, using automatic grading platforms enables to have more accurate grades for homework. There is weekly homework for this class, consisting of 15-20 problems per week. This would be an equivalent of at least 300 human-hours per week devoted only to grading. In practice, this is hardly the case and graders resort to selecting 1-2 problems to grade. This introduces the problem of becoming bias towards students that happened to have these particular problems correct.

\section{Methodology}\label{sec:met}

We analyze the impact of a highly-structured course on student performance by examining the performance gap in student cohorts of the same course in two different years. These courses were held during the Fall quarter with similar incoming frosh students. 

For this, we study the behavior of students in two precalculus courses, one during the Fall quarter 2018, and the other during the Fall quarter 2019. These courses had a similar high-structure and were taught by the same instructor. Both of these classes were large-enrollment (203 students for the 2018 course and 179 students for the 2019 course). 

We analyze the impact of structure on student performance by studying its development from the beginning until the end of the quarter. This effectively makes the initial performance of each cohort to be their own control group.

The performance gap in students is analyzed by considering the behavior of the standard deviation on the normalized cumulative student grades throughout the quarter. In this way, we analyze the behavior of the performance of students as a group, as well as the impact of structure in lowering the difference in experiences that students could be exposed to. 

There were a total of 94 assignments for the 2018 class and 96 for the 2019. These include daily in-class activities, online reflections, weekly worksheets, discussion section quizzes, exams, and projects, among others. The number of assignments made it possible to have this time-dependent analysis.

The time series analyzed is composed by the standard deviation of the normalized cumulative grades up to certain assignment of the entire class. For this, let $X_n$ be the distribution of grades -in percentages- of the $n$th assignment for the class in chronological order. We consider the distribution of the normalized cumulative grades by considering, 
$$\hat{X}_n=\left(\frac{1}{\sum_{k=1}^nw_k}\right)\sum_{k=1}^n w_kX_k\,,$$
where the $k$th assignment is worth $w_k$ percent of the total grade in the course. Let $s_n$ be the standard deviation of $\hat{X}_n$. We are interested in finding a trend of the $s_n$ to determine if there is an impact on the performance gap throughout the quarter. For this we analyze the behavior of $s_n$ by performing a trend analysis using the Mann-Kendall Trend Test \cite{kendall1938new,mcleod2005kendall,mann1945nonparametric,Hussain2019}. 

\section{Results}\label{sec:res}

\begin{figure}[h]
  \centering
  \begin{subfigure}[b]{0.45\linewidth}
    \includegraphics[width=\linewidth]{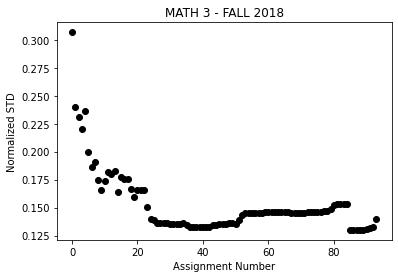}
    \caption{2018 data}
  \end{subfigure}
  \begin{subfigure}[b]{0.45\linewidth}
    \includegraphics[width=\linewidth]{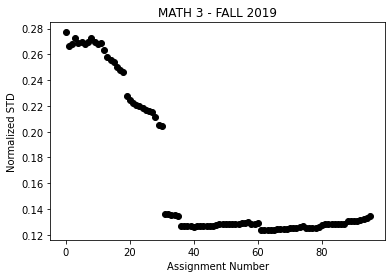}
    \caption{2019 data}
  \end{subfigure}
  \caption{Standard deviation of normalized cumulative grades for each course throughout the quarter. }
  \label{fig::tseries}
\end{figure}

The time series for the standard deviation of the normalized cumulative grades for both groups are shown in \Cref{fig::tseries}. It is important to remark that these time series are labeled by assignment number in chronological order and not by date. Since the courses had at least one assignment per day, the analysis of these time series accurately reflects the development of student performance throughout the quarter.

The results of the trend analysis can be found in \Cref{tab:trend}. For both courses, the Mann-Kendall trend test rejects the null hypothesis that the data does not present any trends with the alternative hypothesis of a decreasing trend. 

These results show that the high-structure used in both courses impacted the grade distribution, considerably reducing the initial standard deviation gap. This implies that towards the end of the quarter, the performance of students became less disperse. 

\begin{table}[h]
\centering
    \begin{tabular}{c c c c c c}
    \toprule
        \multirow{2}{*}{Course} & \multicolumn{5}{c}{Mann-Kendall Trend Test}\\
               &  S & z-value & $\tau$-value & $p$-value & Trend\\
         & \multicolumn{5}{c}{\rule[0.2cm]{260pt}{0.5pt}} \\
        MATH 3 - FALL 2018 & -1105.0 & -3.606 & -0.253 & 0.00031 & Decreasing\\
        MATH 3 - FALL 2019 & -2112.0   & -6.682 & -0.463 & $2.360\times 10^{-11}$ &    Decreasing \\
        \bottomrule
    \end{tabular}
    \caption{Trend Analysis}
    \label{tab:trend}
\end{table}

\section{Conclusions}\label{sec:con}
Large-enrollment courses present several intrinsic challenges for both instructors and students. Effective teaching and learning become difficult as the number of students increase. 

A good framework to assess and address these issues is the use of high-structure on the design and implementation of large-enrollment courses. Highly-structured courses include formative assessments, active learning activities, and low-stakes assignments. These practices can be achieved even in large-enrollment courses with the use of technology and teaching teams. 

Moving into a high-structure region can improve the performance in a course of students as a whole. Addressing scalability, engagement, and feedback in a course with low-stakes assignments, active learning activities, and formative assessments, lowers the performance gap in large-enrollment courses.

\bibliographystyle{plain}
\bibliography{bibliography.bib}

\end{document}